\documentclass[11pt,a4paper]{article}

\usepackage[T1]{fontenc}
\usepackage[utf8]{inputenc}
\usepackage{authblk}
\usepackage{vmargin}
\usepackage{amsthm, amsmath, amssymb}
\usepackage{natbib}
\usepackage{times}
\usepackage{bm}
\usepackage[pdftex]{graphicx}
\usepackage[pdftex]{color}
\usepackage[normalem]{ulem}
\usepackage{paralist}
\usepackage{subfigure}
\usepackage[plain,noend]{algorithm2e}
\usepackage{xr}
\usepackage{dsfont}
\usepackage[final, pdftex]{pdfpages}

\allowdisplaybreaks

\theoremstyle{plain}
\newtheorem{theorem}{Theorem}
\newtheorem{proposition}{Proposition}

\theoremstyle{definition}
\newtheorem{definition}{Definition}

\newtheorem{example}{Example}[section]

\DeclareMathOperator*{\argmin}{\arg\!\min}

\DeclareMathOperator*{\cov}{\rm cov}
\DeclareMathOperator*{\tr}{\rm tr}

\DeclareMathOperator*{\E}{\rm \mathbb{E}}

\newcommand{\muh}{\hat{\mu}}

\newcommand{\srcomment}[1]{}

\newcommand{\srdel}[1]{}
\newcommand{\srdelmath}[1]{}
\newcommand{\srdelfig}[1]{}

\newcommand{\skcomment}[1]{}

\newcommand{\skdel}[1]{}
\newcommand{\skdelmath}[1]{}
\newcommand{\skdelfig}[1]{}

\title{When Does More Regularization Imply Fewer Degrees of Freedom? Sufficient Conditions and Counter Examples from Lasso and Ridge Regression}

\author[1]{Shachar Kaufman\thanks{shachark@post.tau.ac.il}}
\author[1]{Saharon Rosset\thanks{saharon@post.tau.ac.il}}
\affil[1]{Department of Statistics and Operations Research, Tel Aviv University}

\begin{document}

\maketitle

\begin{abstract}
Regularization aims to improve prediction performance of a given statistical
modeling approach by moving to a second approach which achieves worse training
error but is expected to have fewer degrees of freedom, i.e., better agreement
between training and prediction error. We show here, however, that this expected
behavior does not hold in general. In fact, counter examples are given that show
regularization can increase the degrees of freedom in simple situations,
including lasso and ridge regression, which are the most common regularization
approaches in use. In such situations, the regularization increases both training
error and degrees of freedom, and is thus inherently without merit. On the other 
hand, two important regularization scenarios are described where the expected
reduction in degrees of freedom is indeed guaranteed: (a) all symmetric linear
smoothers, and (b) linear regression versus convex constrained linear regression
(as in the constrained variant of ridge regression and lasso).
\end{abstract}

\section{Introduction}
\label{sec:intro}

Let $y \in \mathds{R}^n$ be a random data set generated according to a
probability distribution $f(y;\mu)$, where $\mu = \E(y)$ is a parameter that we
wish to model. A modeling approach $\mathcal{M}: \mathds{R}^n \mapsto
\mathds{R}^n$ is a mapping from a training set $y$ to a model ${\muh}$. Models
are evaluated according to an error or loss criterion $\mathcal{E}(\muh,
y^{\mbox{new}})$, where $y^{\mbox{new}}$ (the test set) is also drawn from $f$,
independently of $y$. Here we focus on the squared-error criterion:
\begin{equation}\label{eq:squared_loss}
   \mathcal{E}(\muh,y^{\mbox{new}}) = \frac{1}{n} \sum\limits_{i=1}^n{\bigl(\muh_i - y^{\mbox{new}}_i\bigr)^2},
\end{equation}
although as will be discussed in Section~\ref{sec:conclusion}, other choices are
possible, with most of the theory that follows intact. We follow previous work
in examining the in-sample error, where any covariate values are the same for the training and testing data \citep{Efron1983, Hastie2009}.

In the model selection problem we are given a collection of candidate modeling
approaches, and our goal is to select an approach having small risk, i.e., small
expected prediction error. A typical setting in which this problem arises is
regularization, where a family of {\em nested} modeling approaches is considered
\citep{Hastie2009}. To apply a modeling approach from the family, one first has
to specify the value of a tuning parameter which controls the amount of fitting
to training data. Consider for example the problem of estimating $\mu$ by
fitting a polynomial to a set of observations $\{(y_i,x_i)\}_{i=1}^n$ with least
squares. Here the degree of the polynomial, $p$, plays the role of the tuning
parameter, where a higher degree leads to more fitting to training data than
does a lower degree. Having specified $p$, the least squares optimization
problem over the set of polynomials of this degree constitutes a modeling
approach, and the solution given specific training data is a model. Choosing a
very low degree (underfitting) is undesirable, as some of the information that
could be gained from the training data is wasted. A high degree often leads to
high variance and overfitting, and is also undesirable. Both underfitting and
overfitting lead to models achieving suboptimal risk.

Model selection is facilitated by producing estimates for the risk of each
candidate modeling approach
\citep{Mallows1973,Akaike1974,Schwarz1978,Stone1974}. The training error,
$\mathcal{E}(\muh,y)$, is a naive such estimate which is typically negatively
biased due to the fact that fitting has been carried out on the same data used
for performance evaluation. This bias (i.e., the difference between the expected
training error and risk) is termed expected {\em optimism} in \citet{Efron1983},
or (up to a constant) {\em effective degrees of freedom} in \citet{Hastie2009}.
Because within a nested family of regularized models the training error
increases monotonically with the level of regularization, the latter is premised
on the idea that the effective degrees of freedom must correspondingly decrease
(just like the degrees of freedom of nested linear regression models decrease
when covariates are removed). In other words, for regularization to have the
potential to reduce risk, less-fitting models (with larger training error) must
have smaller optimism.

We prove in Section~\ref{sec:sufficient} that this is indeed true in some
important cases. However, our main observation in this paper is that such
monotonicity does not necessarily hold. In particular, the lasso
\citep{Tibshirani1996} and ridge regression \citep{Hoerl1962} are two approaches
of great practical and theoretical importance in regularized modeling. As we
show in Section~\ref{sec:examples}, both of them admit counter examples where
more regularized, nested models, also have higher optimism. In fact,
specifically for the lasso this can be argued to be a typical case that arises
in natural examples. Thus monotonicity of degrees of freedom of nested models
breaks down in perhaps the most common and important cases. When the
monotonicity does not hold, the implication is that adding regularization
counter-intuitively increases both the training error and the optimism, and
hence is inherently without merit.

The remainder of the paper is organized as follows. Section~\ref{sec:concepts}
gives formal definitions for the notion of nesting, and reviews the basic 
concepts of optimism, degrees of freedom and their statistical properties.
Section~\ref{sec:sufficient} gives sufficient conditions under which the
effective degrees of freedom grow monotonically in the direction of nesting.
Section~\ref{sec:examples} gives realistic examples where familiar nested models
exhibit reverse monotonicity. Section~\ref{sec:conclusion} offers a
discussion of our results.

\section{Nesting and optimism}
\label{sec:concepts}

\subsection{Nesting}
\label{subsec:nesting}
In traditional linear regression methodology, least squares modeling approaches
are projections on linear subspaces, and nested models are naturally defined
according to the nesting structure of these subspaces. Specifically, one linear
modeling approach is nested in another if the former admits a subset of the
explanatory variables used in the latter. The span of this reduced explanatory
set (defined as the linear span of the observed covariate vectors in 
$\mathds{R}^n$) is geometrically nested in the span of the larger set. Here
we formalize and generalize this definition to cover typical regularized
modeling settings. Our first definition of {\em strict-sense nesting} is an
immediate generalization of the least squares projection notion above.
Intuitively, one modeling approach is nested in another if both fit the model by
minimizing a loss criterion on the training data, over geometrically nested sets
of candidate models (as done in empirical risk minimization~\citep{Vapnik2000}).

\begin{definition}[Strict-sense nesting]
\label{def:strict_sense_nesting} Let $\mathcal{M}_S$ and
$\mathcal{M}_L$ be two modeling approaches that produce models by
optimizing the training error. $\mathcal{M}_S$ performs this
optimization over the {\em model set} $S \subseteq \mathds{R}^n$
while $\mathcal{M}_L$ considers the model set $L \subseteq
\mathds{R}^n$:
\begin{gather*}
  \begin{array}{ccc}
    \mathcal{M}_S: \hat{\mu}^{y,S} = \underset{\tilde{\mu} \in S}{\argmin}\{ \mathcal{E}(\tilde{\mu}, y) \}, & \quad &
    \mathcal{M}_L: \hat{\mu}^{y,L} = \underset{\tilde{\mu} \in L}{\argmin}\{ \mathcal{E}(\tilde{\mu}, y) \}
  \end{array}
\end{gather*}

We say that $\mathcal{M}_S$ is nested in the strict sense in $\mathcal{M}_L$ if $S \subseteq L$. In this case we write $\mathcal{M}_S \preceq \mathcal{M}_L$.
\end{definition}

As elaborated below, strict-sense nesting covers some interesting
regularization families, but others (like penalized ridge regression
or lasso) are not covered by this definition, since different
regularization levels modify the loss criterion rather than the set
of candidate models. Towards this end we devise a second looser
definition, which we term {\em wide-sense nesting}. To be nested in
this sense, the two modeling approaches have to be equivalent to
strict-sense nested approaches for every specific training set $y$, but this
correspondence can be data-dependent.

\begin{definition}[Wide-sense nesting]
\label{def:wide_sense_nesting} Let $\mathcal{Q}$ be a sequence of
nested sets. Let $\mathcal{M}_S$ and $\mathcal{M}_L$ be modeling
approaches that, given a training set $y$, produce the models
$\hat{\mu}^{y,S}$ and $\hat{\mu}^{y,L}$, respectively.
We say that $\mathcal{M}_S$ is nested in the wide sense in
$\mathcal{M}_L$ if, for every value of $y$, there exist sets that
depend on $y$, $Q_S(y) \in \mathcal{Q}$ and $Q_L(y) \in
\mathcal{Q}$, such that $Q_S(y) \subseteq Q_L(y)$ and such that the
models $\hat{\mu}^{y,S}$ and $\hat{\mu}^{y,L}$ are equivalent to the
result of optimizing the same criterion $\mathcal{E}(\tilde{\mu},
y)$ over the model sets $Q_S(y)$ and $Q_L(y)$, respectively. Thus,
the following holds for every value of $y$:
\begin{gather*}
  \begin{array}{ccccc}
      \hat{\mu}^{y,S} = \underset{\tilde{\mu} \in Q_S(y)}{\argmin}\{ \mathcal{E}(\tilde{\mu}, y) \}, & \quad &
      \hat{\mu}^{y,L} = \underset{\tilde{\mu} \in Q_L(y)}{\argmin}\{ \mathcal{E}(\tilde{\mu}, y) \}, & \quad &
      Q_S(y) \subseteq Q_L(y);
  \end{array}
\end{gather*}
In this case it is said that $\mathcal{Q}$ induces the nesting.
\end{definition}

Clearly, by taking $\mathcal{Q}$ to contain the sets $S, L$
themselves, Definition~\ref{def:wide_sense_nesting} is a
generalization of Definition~\ref{def:strict_sense_nesting}.
The two definitions are embodied in the following example:

\begin{example}[ridge regression]
\label{ex:ridge_regression}
The ridge regression modeling approach \citep{Hoerl1962}, in its
common penalized form, fits a model by optimizing a criterion
that incorporates a penalty term weighted by a prespecified tuning parameter
$\lambda \geq 0$:
\begin{align}\label{eq:ridge}
  && \hat{\mu} = \underset{\tilde{\mu} \in S}{\argmin}\Bigl\{\|y - \tilde{\mu}\|_2^2 + \lambda\|\tilde{\beta}\|_2^2 \Bigr\},
  && S = \bigl\{ \tilde{\mu} \mid \exists \tilde{\beta} \in \mathds{R}^p : \tilde{\mu} = X\tilde{\beta}
  \bigr\}.
\end{align}

The Lagrangian dual problem \citep{Boyd2009} of \eqref{eq:ridge} is the less common but
conceptually important constrained form ridge regression. The
dual tuning parameter is $s > 0$, which this time directly constrains
the $L_2$ norm squared of the coefficient vector:
\begin{align}\label{eq:ridge_constrained}
  &&\hat{\mu} = \underset{\tilde{\mu} \in S}{\argmin}\Bigl\{ \|y - \tilde{\mu}\|_2^2 \Bigr\},
  && S = \bigl\{ \tilde{\mu} \mid \exists \tilde{\beta} \in \mathds{R}^p, \|\tilde{\beta}\|_2^2 \leq s : \tilde{\mu} = X\tilde{\beta} \bigr\}.
\end{align}
For the ridge regression problem, the duality of the two forms essentially means
that for a given vector $y$, for each value of $\lambda \geq 0$ there exists a
value of $s > 0$ such that the two problems are equivalent, i.e., give the same
$\muh$ \citep[for details see][]{Davidov2006}.

In the constrained form, the constraint on $\beta$ defines a $p$-ball where all
coefficient vectors must lie. Consequently the model set which is defined by the
projection of this $p$-ball by the matrix $X$ is enclosed by a hyper-ellipsoid
(embedded in the hyperplane spanned by $X$) that scales isotropically with the
value of the tuning parameter $s$. Specifically, according to
Def.~\ref{def:strict_sense_nesting}, a constrained ridge regression with smaller
$s$ is nested in a constrained ridge regression with a larger $s$.

In the penalized form nesting exists in the sense of
Def.~\ref{def:wide_sense_nesting}. Consider two penalized ridge
cases: one with $\lambda_S$ and the other with $\lambda_L$, and
assume $\lambda_S > \lambda_L$. The criterion optimized in
\eqref{eq:ridge} depends on the value of $\lambda$, so there is no
strict sense nesting as per Def.~\ref{def:strict_sense_nesting}.
Now, let
\begin{equation*}
  \mathcal{Q} = \Bigl\{ Q_s = \bigl\{ \tilde{\mu} \mid \exists
  \tilde{\beta} \in \mathds{R}^p, \|\beta\|_2^2 \leq s : \tilde{\mu} =
  X\tilde{\beta} \bigr\} \mid s > 0 \Bigr\}.
\end{equation*}
From the strong duality with the constrained form, the penalized ridge
regression for a given $y$ and $\lambda$ in fact optimizes the
squared-loss criterion over model sets in $\mathcal{Q}$.
In our case, denote these by $Q_S(y)$ and $Q_L(y)$. Since $\lambda_S
\geq \lambda_L$ in the dual (penalized) form, then for every value
of $y$, $\|\hat{\beta}_S\|_2^2 \leq \|\hat{\beta}_L\|_2^2$ and thus
$s_S(y) \leq s_L(y)$ in the primal (constrained) form. Therefore
$Q_S(y) \subseteq Q_L(y)$, and penalized ridge regression models are
nested in the wide sense of Def.~\ref{def:wide_sense_nesting}.
\end{example}

\subsection{Optimism and effective degrees of freedom}
\label{subsec:optimism}

Let the expected training error (averaged over all random training sets) be
$E_{train}$, and the risk or expected prediction error over all
training and test sets be $E_{pred}$. The expected
optimism~\citep{Efron1983}, $\omega$, is then defined as the
difference of these two quantities:
\begin{equation}\label{eq:expected_optimism}
    \omega = E_{pred} - E_{train} =  \frac{1}{n} \left( \E(\|\hat{\mu}- y^{\mbox{new}}\|^2_2) - \E(\|\hat{\mu}-y\|^2_2)\right).
\end{equation}
The optimism theorem, due to \citet{Efron1983, Efron2004}, relates optimism to
the self-influence of observations:
\begin{equation}\label{eq:optimism_theorem}
    \omega = \frac{2}{n} \sum\limits_{i=1}^n{\cov\Bigl(\hat{\mu}_i, y_i\Bigr)}.
\end{equation}
Stein's Lemma~\citep{Stein1981} further states that, under certain regularity
conditions and for a normally distributed, homoscedastic and uncorrelated data set
$y\sim\mathcal{N}(\mu, \sigma^2I)$ (with extensions for a variety of other cases as well  \citep{kumar2009stein}), the expected optimism is proportional to the
divergence (the trace of the Jacobian) of $\hat{\mu}$ as a function of $y$:
\begin{equation}
  \label{eq:steins_lemma}
  \omega = \frac{2\sigma^2}{n} \E\Bigl(\sum\limits_{i=1}^n{\frac{\partial \hat{\mu}_i}{\partial y_i}}\Bigr).
\end{equation}
Expected optimism can thus be thought of as a sensitivity measure of the
fitted values to their respective observation. A short review of the importance
of optimism for model selection is given in
Section~\ref{sec:additional_background} of the Supplementary Material.

Specifically, for the class of linear smoothers of the form
$\hat{\mu} = Sy$, where $S$ is uncorrelated with $y$ (or is simply
fixed), if we assume $\cov(y) = \sigma^2I$,
then~\eqref{eq:optimism_theorem} allows us to directly derive the
optimism as:
\begin{equation*}
  \label{eq:optimism_of_a_linear_smoother}
  \omega = \frac{2\sigma^2}{n}\tr(S).
\end{equation*}
In the case of linear regression and penalized ridge regression, the
matrix $S$ has the form:
\begin{equation}
  \label{eq:ridge_solution}
  \muh = X(X^{T}X + {\lambda}I_{p\times p})^{-1}X^{T}y,
\end{equation}
and the optimism of these approaches can be expressed using the singular value
decomposition (SVD) of the design matrix $X = VDU^{T}$:
\begin{equation}
  \label{eq:optimism_of_ridge}
  \omega = \frac{2\sigma^2}{n} \sum_{j=1}^p{\frac{d_j^2}{d_j^2+\lambda}}
\end{equation}
\citep[for the details, see][]{Hastie2009}. In linear regression ($\lambda = 0$)
this simplifies to $\omega = 2\sigma^2p/n$, thus the optimism here is
proportional to the degrees of freedom, which are the number of optimized
parameters in the linear model. This means that for nested linear regressions,
adding explanatory variables indeed increases optimism. Similar monotonicity
occurs in penalized ridge regression with a general $\lambda > 0$ since~\eqref{eq:optimism_of_ridge} decreases as $\lambda$ increases.

These results motivated the definition by \citet{Wahba1983} of $\tr(S)$ as
``equivalent degrees of freedom'' for linear smoothers. A natural extension is
the definition of ``generalized degrees of freedom'' in \citet{Ye1998} or
``effective degrees of freedom'' in \citet{Hastie2009} for an arbitrary modeling
approach based on the concept of expected optimism:
\begin{equation*}
  \label{eq:gdf}
  \mathrm{df} = \omega\frac{n}{2\sigma^2}
  \overset{(*)}{=} \frac{1}{\sigma^2} \sum\limits_{i=1}^n{\cov\bigl(\hat{\mu}_i, y_i\bigr)}
  \overset{(**)}{=} \E\Bigl(\sum\limits_{i=1}^n{\frac{\partial \hat{\mu}_i}{\partial y_i}}\Bigr).
\end{equation*}
When applicable, the equality $(*)$ follows from the optimism
theorem~\eqref{eq:optimism_theorem}, while $(**)$ comes from Stein's
lemma~\eqref{eq:steins_lemma}. In the original linear regression context, the
degrees of freedom are a measure of both the optimism, and of the amount of
regularization (implying for any pair of models which model is nested in the
other). As shown above for penalized ridge regression, and as will be shown more
generally in the theorems of the next section, a monotonic nondecreasing
relation between the amount of regularization and the effective degrees of
freedom also holds in other important regularization methods. This belies a
notion that such monotonicity holds in general, providing a wide theoretical
basis for applying regularization. However, this is not always the case.

Before going into the examples of Section~\ref{sec:examples} which
deal with lasso and (constrained) ridge regression and are of more practical
relevance, let us begin with a simple illustrative counter-example where
strict-sense nesting does not imply monotonicity in optimism.

\begin{example}[toy counterexample]
\label{ex:line_inside_circle}

Assume our data vector is two dimensional: $y\in \mathds{R}^2$.  Let $S$ be a
line segment in $\mathds{R}^2$ where its first coordinate is in $[-1,1]$ and its
second is $0$. Let $L$ be the unit disk. We shall relax this later, but first
let our data be $y_1 \sim \mathcal{U}(-1,1)$, and $y_2=2$. Projecting from any
realization of $y$ on the line segment $S$ gives $\muh^S = (y_1, 0)$. On the
other hand, projecting onto the disk $L$, the correlation in the first
coordinate is partial (and since $y_2$ is fixed, correlation of the second
component is still zero). Formally, the optimism for these models can be
computed as:

\begin{align*}
    y &= \begin{bmatrix} 0 \\ 2 \end{bmatrix} + \begin{bmatrix} \varepsilon \\ 0 \end{bmatrix} &\varepsilon &\sim \mathcal{U}[-1,1] \\
    \muh^S &= \left\{\begin{array}{cl}
      \begin{bmatrix} \makebox[20pt]{$ -1$} \\ 0 \end{bmatrix}, & \text{if $y_1 < -1$} \\[+10pt]
      \begin{bmatrix} \makebox[20pt]{$  1$} \\ 0 \end{bmatrix}, & \text{if $y_1 > 1$} \\[+10pt]
      \begin{bmatrix} \makebox[20pt]{$y_1$} \\ 0 \end{bmatrix}, & \text{otherwise}
    \end{array}\right. &
    \muh^L &= \left\{\begin{array}{cl}
      y, & \text{if $\|y\| < 1$} \\[+5pt]
      \frac{y}{\|y\|}, & \text{otherwise}
    \end{array}\right.\;,
\end{align*}
where only the latter cases of each model are relevant for our distribution.


\begin{align*}
    \cov(\muh_1^S, y_1) &= \E[(\muh_1^S-0)(y_1-0)] = \E(y_1^2) = \frac{1}{3} \\
    \cov(\muh_2^S, y_2) &= \cov(\muh_2^L, y_2) =0 \\
    \cov(\muh_1^L, y_1) &= \E\Bigr[(\frac{y_1}{\sqrt{y_1^2+4}}-0)(y_1-0)\Bigl]  \\
    &= \frac{1}{2} \Bigr(\frac{\sqrt{5}}{2} - 2 \log{(1+\sqrt{5})} - \frac{-\sqrt{5}}{2} + 2 \log{(-1+\sqrt{5})}\Bigl) \approx 0.1556.\\
    \omega^S &= \frac{2}{n} \sum_{i=1}^n{ \cov(\muh_i^S, y_i) } = \frac{1}{3}\;,\;\;
    \omega^L = \frac{2}{n} \sum_{i=1}^n{ \cov(\muh_i^L, y_i) } \approx 0.1556,
\end{align*}
and indeed, the smaller nested model set, $S$, leads to more optimism than the larger one, $L$.

Supplementary Figure~\ref{fig:ex_line_ball} (right panel) demonstrates the
phenomenon using a Monte-Carlo simulation with $10^6$ draws of $y$. When $y$ is
normally distributed having the same expectation and variance-covariance matrix
we still see that the larger approach has smaller optimism. This phenomenon is
not limited to a two-dimensional setup; supplementary
Figure~\ref{fig:ex_line_ball_highdim} shows its persistence when the disk is
replaced by a $n$-ball and the line segment is replaced by a hyperplane tile.
Furthermore, Supplementary Section~\ref{sec:toy} shows that such nonmonotonicity
in optimism can have a significant effect on prediction error.
\end{example}

\section{Sufficient Conditions}
\label{sec:sufficient}

We propose in this section two theorems which address important special cases of
nesting. We show first (Thm.~\ref{thm:linear_smoothers}) that for the class of
symmetric linear smoothers (spanning most commonly used smoothing approaches),
nesting in the wide sense (Def.~\ref{def:wide_sense_nesting}) guarantees the
smaller modeling approach has less optimism. Our next result
(Thm.~\ref{thm:sufficient_conditions}) concerns the case where the smaller
modeling approach is a projection on a convex set $S$, while the bigger one is a
projection on a linear subspace containing $S$. Here the modeling approaches are
nested in the strict sense (Def.~\ref{def:strict_sense_nesting}), and we show
that in this case as well monotonicity of optimism is guaranteed.

\begin{theorem}\label{thm:linear_smoothers}
Let $y$ be a homoscedastic and mutually uncorrelated observation
vector. Let $\muh^S = Sy$, $\muh^L = Ly$ be two
linear smoothers that are real and symmetric ($S^T=S$, $L^T=L$). If
$\muh^S$ is nested in $\muh^L$ in the wide-sense over convex sets as
in Definition~\ref{def:wide_sense_nesting}, then $\muh^L$ has more
optimism than $\muh^S$.
\end{theorem}

A proof is provided in Appendix~\ref{apdx:linear_smoothers_proof}. The 
following example uses Thm.~\ref{thm:linear_smoothers} to show 
that all nested generalized ridge regressions exhibit monotone optimism
in the direction of nesting.

\begin{example}[generalized ridge regression]
Consider the following family of modeling approaches:
\begin{gather}\label{eq:genridge}
  \hat{\mu} = X \underset{\beta \in \mathds{R}^p}{\argmin}\Bigl\{\|y - X\beta\|_2^2 + \lambda\beta^T K \beta \Bigr\},
\end{gather}
where $K$ is some symmetric matrix. The solution is
$\hat{\mu}(\lambda) = X(X^TX + \lambda K)^{-1} X^Ty$. It can easily
be verified that $\hat{\mu}(\lambda)^T XKX^T\hat{\mu}(\lambda)$ is a
monotone decreasing function of $\lambda$. Thus we have nesting in
the wide sense over the sets
\begin{equation*}
  \mathcal{Q} = \Bigl\{ Q_s = \bigl\{ \tilde{\mu}:
  \tilde{\mu}^TXKX^T\tilde{\mu}\leq s\bigr\} \mid s > 0 \Bigr\},
\end{equation*}
and according to Thm. \ref{thm:linear_smoothers} the optimism is
indeed monotone in $\lambda$.

Although direct eigen-analysis can give an explicit derivation of the optimism
and therefore also prove the monotonicity in special cases (including ridge
regression as shown in Section~\ref{sec:concepts}, and natural smoothing splines
as in \citet{Hastie2009}), to our knowledge there is no previous general result
that can be used to prove monotonicity for all generalized ridge approaches.
\end{example}

\begin{theorem}
\label{thm:sufficient_conditions}
Let $\mathcal{M}_S \preceq \mathcal{M}_L$ be nested modeling
approaches as defined by Definition~\ref{def:strict_sense_nesting}
with the squared error criterion
\eqref{eq:squared_loss}. Let $L$ be a linear subspace of
$\mathds{R}^n$ and $S \subseteq L$ a convex set. If the
conditions for Stein's lemma \eqref{eq:steins_lemma} are satisfied
for both $\mathcal{M}_S$ and $\mathcal{M}_L$, then $\omega_S \leq
\omega_L$.
\end{theorem}

This theorem, proven in Appendix~\ref{apdx:sufficient_conditions_proof}, implies
that any constrained linear regression model including constrained ridge
regression, constrained lasso (discussed below), constrained elastic
net~\citep{zou2005regularization} and others, has lower optimism and fewer
degrees of freedom than the unconstrained linear regression model with the same
variables.

\begin{example}[convexity requirement in Thm.~\ref{thm:sufficient_conditions}]
\label{ex:convexity_requirement}
It is important to note that the convexity requirement is needed in
Thm.~\ref{thm:sufficient_conditions}. Let $$ y = \begin{bmatrix} 1 \\ 0
\end{bmatrix} + \varepsilon\;,\; \varepsilon \sim
\mathcal{N}\left(\begin{bmatrix} 0 \\ 0 \end{bmatrix}, (0.1)^2 I_{2 \times
2}\right),$$ be a two-dimensional data distribution. Let $L$ be the vertical
axis ($y_1 = 0$) and let $S$ be a two-point set $\{(0,-1), (0,1)\}$. Let $\mathcal{M}_S$ and
$\mathcal{M}_L$ be Euclidean projections on $S$ and $L$ respectively. Their fits are
$\hat{\mu}^S_2 = \mbox{sgn}(y_2)$ and $\hat{\mu}^L_2 = y_2$, with the first
coordinate of both $\hat{\mu}^S$ and $\hat{\mu}^L$ equal $0$. $S\subset L$ and
therefore $\mathcal{M}_S \preceq \mathcal{M}_L$. But it is easy to verify
that:
\begin{align*}
  \omega_L &= \mbox{var}(y_2) = 0.01 \\
  \omega_S &= \E|y_2| = 0.1 \sqrt{2/\pi} \approx 0.08.
\end{align*}
Thus, the theorem does not hold without the convexity requirement.
\end{example}

\section{Counterexamples}
\label{sec:examples}

Examples \ref{ex:convexity_requirement} and \ref{ex:line_inside_circle} were
simple illustrations that a smaller nested approach can give higher optimism
than a larger one. While it is easy to devise such anecdotal examples, a key
question is to what extent may we expect to encounter this phenomenon in the
wild, i.e., in practically interesting and relevant situations. To address this
question we call on what are perhaps the two most widely used and studied
regularization approaches in regression: lasso \citep{Tibshirani1996} and ridge
regression \citep{Hoerl1962}. Supplementary Section~\ref{sec:irp} gives a more
exotic third example, using regularized isotonic regression.

\begin{example}[lasso counterexample]
\label{ex:lasso_counter_example}

Using similar notations to the ridge definitions above, the
penalized and constrained formulations of lasso are, respectively:
\begin{align}\label{eq:penlasso}
  && \hat{\mu} &= \underset{\tilde{\mu} \in S}{\argmin}\Bigl\{\|y - \tilde{\mu}\|_2^2 + \lambda\|\tilde{\beta}\|_1 \Bigr\},
  && S = \bigl\{ \tilde{\mu} \mid \exists \tilde{\beta} \in \mathds{R}^p : \tilde{\mu} = X\tilde{\beta}
  \bigr\}, \\
  \label{eq:lasso_constrained}
  && \hat{\mu} &= \underset{\tilde{\mu} \in S}{\argmin}\Bigl\{ \|y - \tilde{\mu}\|_2^2 \Bigr\},
  && S = \bigl\{ \tilde{\mu} \mid \exists \tilde{\beta} \in \mathds{R}^p, \|\tilde{\beta}\|_1 \leq s : \tilde{\mu} = X\tilde{\beta} \bigr\}.
\end{align}
Like in ridge regression, constrained lasso modeling approaches are
nested in the strict sense (Def.~\ref{def:strict_sense_nesting}), and penalized 
lasso modeling approaches are nested in the wide sense
(Def.~\ref{def:wide_sense_nesting}).

If we denote the solution of~\eqref{eq:penlasso} by $\hat{\mu}(\lambda) = X
\hat{\beta}(\lambda)$, it is well known that $\hat{\beta}(\lambda)$ corresponds
to soft variable selection, where $\hat{\beta}(\lambda)_j = 0$ for some of $j=1,
\dots, p$. Following \citet{Zou2007}, for a specific lasso solution, we denote
by ${\cal A} \subseteq \{1, \dots, p\}$ the {\em active set} of variables with
non-zero coefficients, i.e., $\hat{\beta}(\lambda)_j = 0,\;\forall j \notin \cal
A.$ To avoid complex notation, the dependence of $\cal A$ on the penalty
$\lambda$ or constraint $s$ is left implicit. \citet{Zou2007} give the Stein
unbiased estimate of optimism for the penalized lasso formulation: $\hat{\omega}
= |\cal A|,$ while \citet{kato2009degrees} gives a slightly different result for the
constrained version. As the regularization level decreases in the lasso solution
of a specific data set, the number of active variables can decrease and not only
increase \citep{Efron2004LARS, Zou2007}. Hence the Stein unbiased estimate of
optimism is in general not expected to be monotone increasing as regularization
decreases. The important question, however, is to what extent can this behavior
exist in expectation over a distribution. In other words, can the optimism
itself adopt a similar pattern in realistic examples?

We offer a simple example which demonstrates that this is eminently
possible. Consider a regression problem with three covariates,
$n = 1001$ observations and the following characterization:
\begin{align*}
  x_{1i} &= \left\{ \begin{array}{ll}
    \sqrt{\frac{1}{n-1}} & \mbox{ if $i<n$} \\
    0 & \mbox{   if $i=n$}
    \end{array} \right. &
  x_{2i} &= \left\{ \begin{array}{ll}
    \sqrt{\frac{1}{n-1}} & \mbox{ if $i$ odd, $i<n$} \\
    -\sqrt{\frac{1}{n-1}} & \mbox{ if $i$ even, $i<n$} \\
    0 & \mbox{ if $i=n$}
  \end{array} \right. \\
  x_{3i} &= \left\{ \begin{array}{ll}
    \sqrt{\frac{3}{2(n-1)}} & \mbox{ if $i$ odd, $i<n$} \\
    0 & \mbox{ if $i$ even, $i<n$} \\
    0.5 & \mbox{ if $i=n$}
  \end{array} \right. &
  &\begin{array}{ll}
    y_i &= x_{1i} + x_{2i} -0.1x_{3i} + \varepsilon_i \\
    \varepsilon_i &\sim \mathcal{N}(0,0.02) \mbox{ i.i.d.}
  \end{array}
\end{align*}

\begin{figure}[htbp]
  \centering
  \mbox{
    \subfigure{\includegraphics[width=0.45\textwidth]{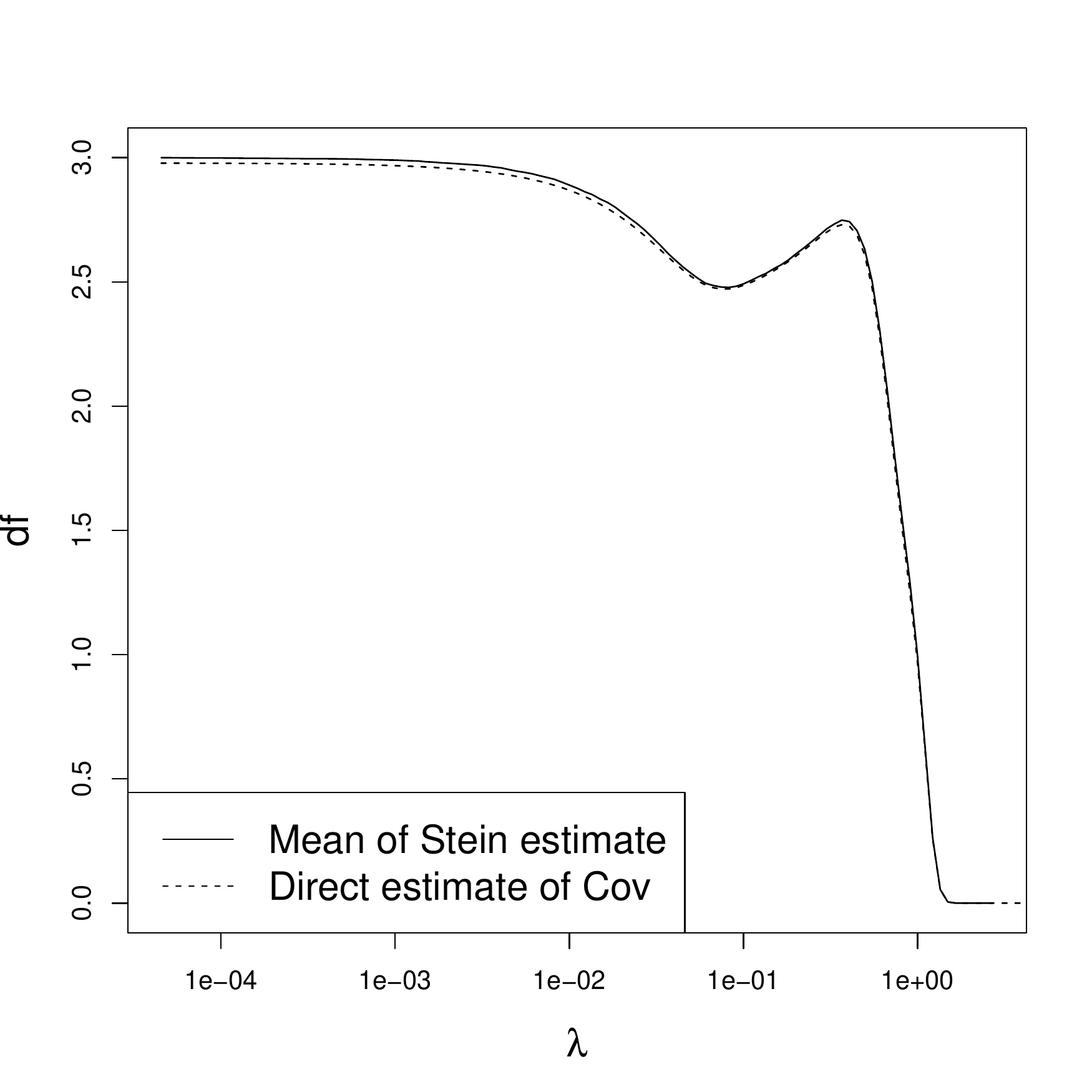}}
    \quad
    \subfigure{\includegraphics[width=0.45\textwidth]{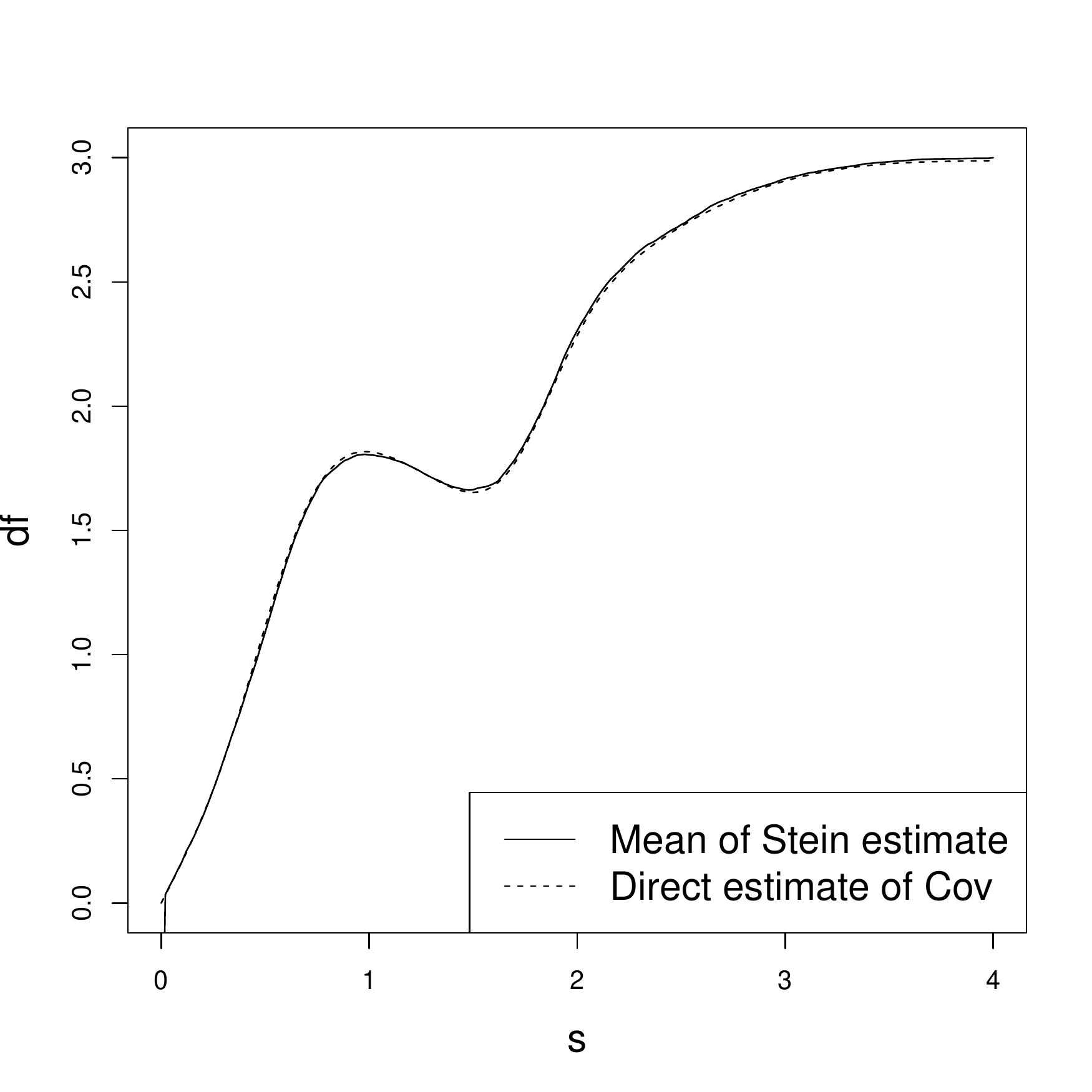}}
  }
  \caption{Lasso counter-example. The plots show the (effective) degrees of freedom
  of the penalized (left) and constrained (right) forms as a
  function of the regularization parameter. The degrees of freedom are calculated
  independently as the empirical expectation of the Stein estimate
  (solid) and the empirical estimate of  $\sum_i \mbox{cov}(y_i,
  \hat{y}_i)$. The two estimates agree and both demonstrate the clear
  non-monotonicity of degrees of freedom in the regularization level. The estimated
  standard errors of the degrees of freedom estimates based on our simulations are
  negligible compared to the non-monotonicity (not shown on the
  graphs). See text for details of the simulation setup.}
  \label{fig:lasso}
\end{figure}

Figure~\ref{fig:lasso} (left panel) shows the degrees of freedom of the
penalized lasso solution for this problem as a function of $\lambda$, estimated
from $5000$ independent simulations in two manners: \begin{inparaenum}[(i)]
\item directly, using the definition of optimism in
Eq.~\eqref{eq:optimism_theorem}, and \item by calculating the Stein estimate
(number of non-zero coefficients) and averaging it over the simulations.
\end{inparaenum} As expected, the two calculations agree well, and both
demonstrate the clear non-monotonicity of the optimism in the regularization
level. The implications for model selection are clear: at $\lambda=0.5$ the
modeling approach has both higher in-sample error and higher optimism than at
$\lambda=0.1$, therefore the value $0.5$ is not a useful value of the
regularization parameter to consider. For this distribution, decreasing the
tuning parameter from $0.5$ to $0.1$ in fact decreases the true measure of
regularization as captured by optimism and degrees of freedom; thus, the
apparently less regularized model is in fact more regularized. If we move from
the penalized formulation to the constrained formulation (using Kato's
derivation of the Stein estimate), the phenomenon persists (Fig.~\ref{fig:lasso}, right panel).
\end{example}

\begin{example}[ridge regression counterexample]
\label{ex:ridge_counter_example}

For penalized ridge regression, the explicit derivation of its
optimism in Eq.~\eqref{eq:optimism_of_ridge} guarantees monotonicity
between the regularization level and the optimism, hence it cannot
admit counter examples. The standard result in
Eq.~\eqref{eq:optimism_of_ridge} assumes homoscedastic error, we now
generalize it to the heteroscedastic case as well:

\begin{proposition}\label{prop:heteroscedastic_optimism}
  Suppose observations $y_i$, which are components of the vector $y_{n \times 1}$, are mutually uncorrelated but not homoscedastic, i.e., the covariance matrix of $y$ has components:
  \begin{equation*}
    \Lambda_{ij} = \begin{cases} \sigma_i^2 &, \text{if $i = j$} \\ 0 &, \text{otherwise} \end{cases}.
  \end{equation*}
  Then the expected optimism of the ridge regression modeling approach is
  \begin{equation*}
    \omega = \frac{2}{n} \sum_{i=1}^n{\sigma_i^2 \sum_{j=1}^p{\frac{d_j^2}{d_j^2 + \lambda} u_{ij}^2}},
  \end{equation*}
  where $d_j$ and $u_{ij}$ are components of the matrices $D$ (on the main diagonal) and $U$ respectively, in the SVD of the design matrix: $X=VDU^T$.
\end{proposition}

We leave the proof to Supplementary Section~\ref{apdx:heteroridge}. \\

On the other hand, such monotonicity does not always hold for the constrained form~\eqref{eq:ridge_constrained}. As mentioned in
Section~\ref{subsec:nesting}, fitting is done by projection onto a model set
that is enclosed by a hyper-ellipsoid centered at the origin. Changing the value
of the regularization parameter $s$ shrinks or inflates the hyper-ellipsoid
isotropically. Surprisingly, in this case too there are setups where we get
smaller optimism when projecting onto a larger ellipsoid. We describe next a
relatively simple setup that demonstrates this.

Let $S$ be the set enclosed by a hyper-ellipsoid centered at the origin, with
principal directions parallel to the axes and equatorial radii $(r_S, r_S,
\dots, r_S, r_Sh)$. Let $L$ be a similarly defined set but with radii $(r_L,
r_L, \dots$ $, r_L, r_Lh)$. The parameter $h < 1$ determines eccentricity in the
last component. This setup is presented in two dimensions in Supplementary
Figure~\ref{fig:ex_ellipse_zoom} (top left panel). Specifically, consider taking
the following parameter values: $n = 2, h = 0.1, r_S = 1, r_L = 10$. These
values can be thought of as constrained ridge regression with a diagonal design
matrix $X$. We begin with an illustrative distribution: $y_2 \sim \mathcal{U}(3,
5), y_1 = 1$. The rationale behind this setup is that $y$ is situated such, that
the image of its Euclidean projection on $L$ is a nearly horizontal line
segment, while its projection on $S$ has a much larger vertical ($y_2$)
component. There is also the contradictory effect of the circumference of the
ellipses at play, but since we have highly eccentric ellipses, it is much less
pronounced. The correlation between observed and fitted values is therefore
higher for $\mathcal{M}_S$ than it is for $\mathcal{M}_L$. Supplementary
Figure~\ref{fig:ex_ellipse_zoom} (bottom left panel) shows realizations of $y$
for $r_L = 10$, and the two right panels depict the corresponding fitted values
on the same scale. The componentwise covariance of $y$ is visibly greater with
${\hat{\mu}}^S$ than it is with ${\hat{\mu}}^L$.

This behavior persists beyond the illustrative setup described so far: it scales
with the dimension $n$ (hyper-ellipsoids, results not shown), and endures if we
take $y$ to be distributed according to a normal distribution that is
uncorrelated but heteroscedastic. Results are difficult to obtain in closed-form
for this case, because, in constrained ridge regression, projections involve the
solution of quartic and higher order polynomial equations. We therefore settle
for an estimate of the optimism, with appropriate confidence intervals.
Figure~\ref{fig:ex_ellipse_zoom_profile} gives the optimism profile for the
normal distribution
  \begin{equation*}
    y \sim \mathcal{N}\Bigl(
    \begin{bmatrix} 3 \\ 10 \end{bmatrix},
    \begin{bmatrix} 0.1 & 0 \\ 0 & 3 \end{bmatrix}
    \Bigr),
  \end{equation*}
and when the larger model set is inflated starting from $r_L = r_S = 1$ up to
$r_L = 10$. This profile is monotonic nondecreasing up to $r_L \approx 2.4$, but
then becomes strictly decreasing for the remainder of the examined range.

\begin{figure}[htbp]
  \centering
  \includegraphics[width=0.6\textwidth]{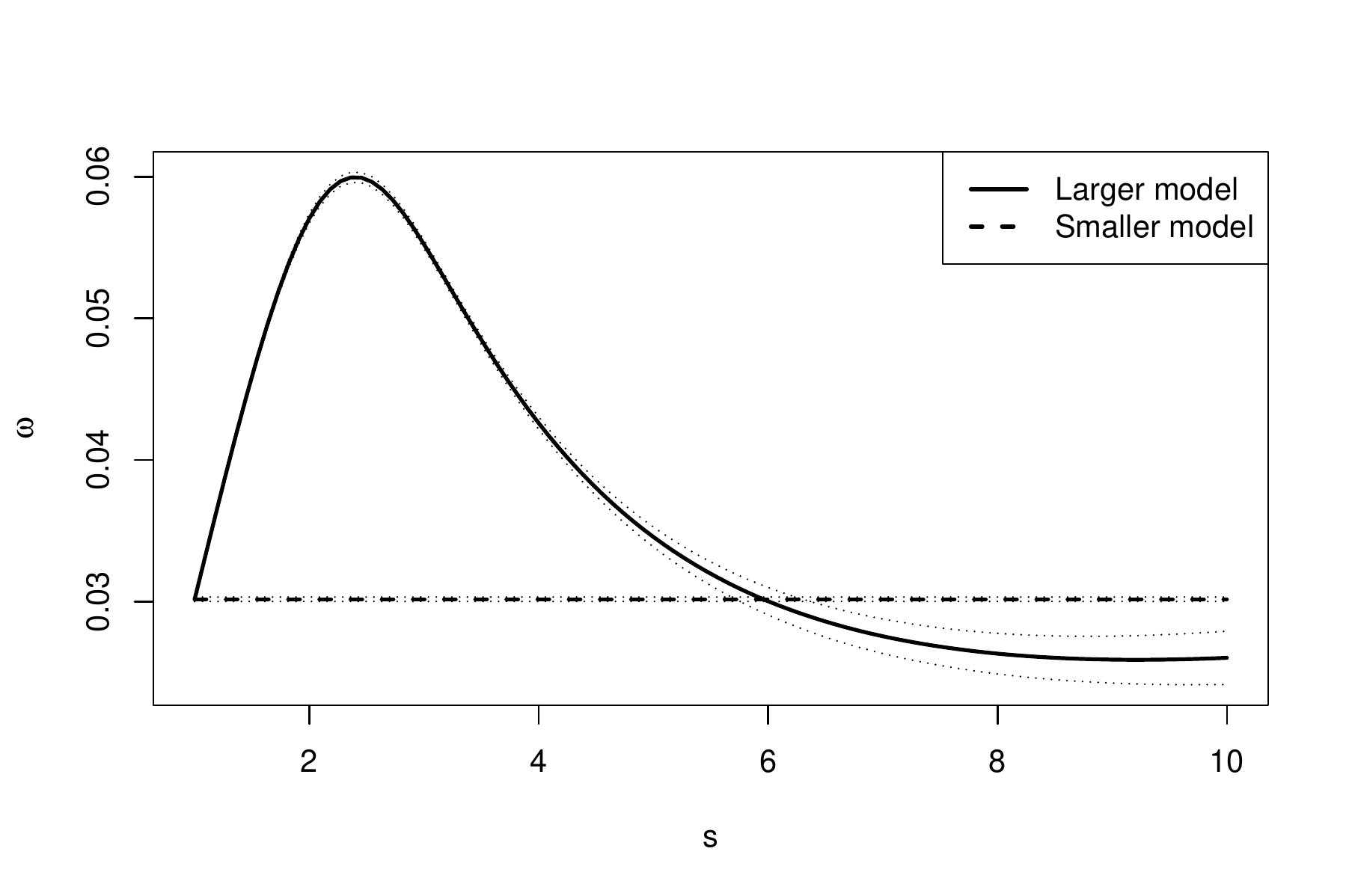}
  \caption{Constrained ridge regression counter-example, when starting
  from $r_S = r_L = 1$ and inflating the larger ellipse until $r_L =
  10$. 95\% confidence intervals for the optimism are shown (dotted).}
  \label{fig:ex_ellipse_zoom_profile}
\end{figure}

Unlike the realistic lasso example, the constrained ridge example above is more
contrived, requiring non-standard error distributions (non-homoscedastic error).
However, the ridge case is potentially more intriguing because the penalized
form guarantees monotonicity, while the constrained form admits counter
examples.
\end{example}

\section{Discussion}
\label{sec:conclusion}

Regularization (nesting) and degrees of freedom (optimism) are key concepts in
statistics and specifically in model selection. Because these concepts are
closely related in the context of fundamental modeling approaches such as
linear regression, a notion of a generally applicable nondecreasing monotonic
relationship between them has permeated the statistics literature
\citep[e.g.][]{Ye1998, Zou2007, Kramer2011}. This notion is also expressed by
the use of optimism to define ``effective degrees of freedom''~\citep{Hastie2009}.
We have shown here that for some important families of nested modeling
approaches, the monotonicity is indeed preserved. However the general
relationship is a misconception that does not hold even in simple and familiar
scenarios such as lasso or constrained ridge regression. In particular, our
lasso (Section~\ref{sec:examples}) and isotonic recursive partitioning
(Section~\ref{sec:irp}) examples are natural and realistic. In such situations,
the fundamental premise of regularization as controlling model complexity and
decreasing optimism is in fact incorrect, and regularized models with more
optimism than their less regularized counterparts are guaranteed to be inferior
in their expected predictive performance.

Specifically for ridge regression with additive heteroscedastic normal noise,
Proposition~\ref{prop:heteroscedastic_optimism} shows that the penalized form
guarantees a monotonic nondecreasing relationship between regularization (as
captured in this case by the tuning parameter, $\lambda$) and optimism,
$\omega$. Surprisingly, the constrained form does not guarantee this, as
demonstrated in Example~\ref{ex:ridge_counter_example}. On the one hand, the two
forms are equivalent (strongly dual) in the optimization theory sense, in that
for every realized training set, we may switch from one form to the other with
an appropriate choice of tuning parameter value, and produce the same model. On
the other hand, this mapping is data-dependent and thus random, which means it
does not imply that the two forms are equivalent statistical modeling
approaches. This subtlety which is reflected in our results also comes up in the
Stein unbiased estimates for penalized and constrained degrees of freedom for the
lasso in \citet{kato2009degrees} and \citet{Zou2007}.

Alternative formulations for the model selection problem exist, which replace
optimism with different notions of complexity, for which monotonicity is
guaranteed. In particular, the machine learning community traditionally defines
model complexity via Vapnik-Chervonenkis (VC) dimension of the model
set~\citep{Vapnik2000}, and calculates penalties on training error which give
bounds on prediction error in place of expected error expressed by optimism
\citep{Cherkassky2007}. The penalties depend monotonically on the VC dimension,
hence the consistency between model complexity and prediction penalty is
guaranteed. A major downside with using this approach is that it gives
worst-case bounds (which are often very loose) in place of estimates of expected
prediction error. More critically, unlike the optimism, these penalties are
independent of both the modeling approach and the true underlying distribution,
depending only on the model set. Hence they are of a fundamentally different
nature than modeling-approach-specific estimates based on optimism.

The concept of optimism is applicable to other loss functions besides the
squared error loss we have focused on here, as shown by \citet{Efron2004}.
Indeed, so is the concept nesting defined in Section~\ref{subsec:nesting}. We
thus expect that the general spirit of our positive results from
Section~\ref{sec:sufficient} and negative results from
Section~\ref{sec:examples} should not change when considering other loss
functions (for example, exponential family log-likelihoods). The details of
these generalizations remain a topic for future research.

\section*{Acknowledgements}
The authors are grateful to E. Aharoni, R. Luss and M. Shahar for
useful ideas and discussion, and to F. Abramovich, T. Hastie, R.
Heller, G. Hooker, R. Tibshirani and the reviewing team for
thoughtful and useful comments. This research was partially
supported by Israeli Science Foundation grant 1487/12 and by a
fellowship to SK from the Edmond J. Safra Center for Bioinformatics
at Tel Aviv University.

\clearpage

\section*{Appendix}

\subsection{Proof of Theorem~\ref{thm:linear_smoothers}}
\label{apdx:linear_smoothers_proof}

\begin{proof}
Let the eigenvectors of $L$ be $\{u_i\}_{i=1}^n$ with associated
eigenvalues $\{\lambda_i\}_{i=1}^n$. Also let the eigenvectors of
$S$ be $\{z_i\}_{i=1}^n$ with associated eigenvalues
$\{\delta_i\}_{i=1}^n$ (some eigenvalues might be zero, some might
not be unique). Since both eigenvector bases are orthonormal and
span $\mathds{R}^n$, we may transform one to the other via a
rotation matrix $R$: $Z=RU$ (where $U$ and $Z$ are matrices whose
rows are individual $u_i$ and $z_i$ respectively).

Since $\muh^S$ is nested in $\muh^L$ in the wide-sense, there exists
a parametrization $\mathcal{Q}$ of $\mathds{R}^n$ with nested
contours, such that for each value of $y$, $Sy$ and $Ly$ are the
Euclidean projection of $y$ onto $Q_S(y) \in \mathcal{Q}$ and
$Q_L(y) \in \mathcal{Q}$ respectively, and such that $Q_S(y)
\subseteq Q_L(y)$. Because of this nesting, $L$ fits a model that is
closer to $y$ than $S$: $\|\muh^S-y\|_2^2 \geq \|\muh^L-y\|_2^2$.
The conditions of this theorem also specify that $Q_L(y)$ and
$Q_S(y)$ are convex.

Since in this case both modeling approaches are linear smoothers,
$\muh^S(y=0)=S\cdot0=0$ and we must have $0 \in Q_S$ (more
generally, the axes origin has to be at the limit of the smallest
contour in the parametrization $\mathcal{Q}$). This also means that
all eigenvalues must be in $[0,1]$, else the projection of their
eigenvector $Sz_i=\delta_i z_i$ cannot be a projection to a convex
set which includes the origin. Since $\|Sz_i-z_i\|_2^2 \geq
\|Lz_i-z_i\|_2^2$, if follows that $Lz_i$ has to be outside the ball
of radius $\delta_i$ around the origin, and hence $\|Lz_i\|_2^2 \geq
\|Sz_i\|_2^2$ (See Figure~\ref{fig:nested_linear_smoothers}).
We therefore have:
  \begin{align*}
    \|Lz_i\|_2^2 &= \Biggl\|L \sum_{j=1}^n{R_{ij} u_j}\Biggr\|_2^2 = \Biggl\|\sum_{j=1}^n{R_{ij} L u_j}\Biggr\|_2^2 =
    \Biggl\|\sum_{j=1}^n{R_{ij} \lambda_j u_j}\Biggr\|_2^2 \\
    &= \sum_{j=1}^n{R_{ij}^2 \lambda_j^2} \geq \|Sz_i\|_2^2 = \delta_i^2.
  \end{align*}
  Subsequently,
  \begin{align*}
    \sum_{i=1}^n{\sum_{j=1}^n{R_{ij}^2} \lambda_j^2} &\geq \sum_{i=1}^n{\delta_i^2} \Longleftrightarrow \\
    \sum_{j=1}^n{\lambda_j^2 \sum_{i=1}^n{R_{ij}^2}} &\geq \sum_{i=1}^n{\delta_i^2}.
  \end{align*}
  But since $R$ is a rotation matrix, the sum of squares along any column or row is unity. Thus,
  \begin{equation*}
    \sum_{j=1}^n{\lambda_j^2} \geq \sum_{i=1}^n{\delta_i^2} \Leftrightarrow \tr(L^T L) \geq \tr(S^T S).
  \end{equation*}
  Let us now reexamine the nesting consequence:
  \begin{align*}
    &\|\muh^S-y\|_2^2 \geq \|\muh^L-y\|_2^2 &\Longleftrightarrow \\
    &y^T (S-I)^T (S-I)y \geq y^T (L-I)^T (L-I)y &\Longleftrightarrow \\
    &(S-I)^T (S-I) - (L-I)^T (L-I) \succeq 0 &\Longleftrightarrow \\
    &(S^T S - L^T L) + (L-S) + (L-S)^T \succeq 0 &\Longrightarrow \\
    &\tr(S^T S) - \tr(L^T L) + 2\tr(L-S) \geq 0.
  \end{align*}
  With the previous result we must have
  \begin{equation*}
    \tr(L) \geq \tr(S).
  \end{equation*}
  Which, for data distributed according to $y \sim \mathcal{F}(\mu,\sigma^2 I$) (i.e., homoscedastic and mutually uncorrelated) and for $L_2$-loss, means that
  \begin{equation*}
    \omega(L) \geq \omega(S).
  \end{equation*}
  Thus, $\muh^L$ has more optimism than $\muh^S$.
\end{proof}

\subsection{Proof of Theorem~\ref{thm:sufficient_conditions}}
\label{apdx:sufficient_conditions_proof}

\begin{proof}
  Jacobian main-diagonal components for the $\mathcal{M}_S$ modeling approach are given by
  \begin{equation*}
    \frac{\partial{\muh_i^{y,S}}}{\partial{y_i}} =
    \lim_{\varepsilon \rightarrow 0}{ \frac{ \muh_i^{y + \varepsilon e_i,S} - \muh_i^{y - \varepsilon e_i,S} }{2\varepsilon} },
  \end{equation*}
  and similarly for $L$ ($e_i$ is the unit vector whose $i$'th component equals 1).

  For every value of $y \in \mathds{R}^n$, for every value of $\varepsilon \in \mathds{R}$ and for each $ i \in \{1,\dots,n\}$, we have
  \begin{equation*}
    \muh_i^{y + \varepsilon e_i,S} - \muh_i^{y - \varepsilon e_i,S} \leq
    {\|\muh_i^{y + \varepsilon e_i,S} - \muh_i^{y - \varepsilon e_i,S}\|}_2.
  \end{equation*}
  A projection mapping onto a convex set is a non-expansion mapping \citep[e.g., as used in][]{Tibshirani2011}. There thus exist $k_S \leq 1$ and $k_L \leq 1$ such that for every two values of $y$: $a$ and $b$,
  \begin{align*}
    {\|\muh^{a,S} - \muh^{b,S}\|}_2 &\leq k_S{\|a - b\|}_2 \\
    {\|\muh^{a,L} - \muh^{b,L}\|}_2 &\leq k_L{\|a - b\|}_2.
  \end{align*}
Because $L$ is a linear subspace, the Euclidean projection onto $S$
may be broken down to first projecting onto $L$ and then projecting
from there onto $S$: $\muh^{y,L,S} \equiv \muh^{y,S}$ (as shown in
Supplementary Figure~\ref{fig:sequential_projection}).
  Hence
  \begin{equation*}
    {\|\muh^{a,S} - \muh^{b,S}\|}_2 = {\|\muh^{a,L,S} - \muh^{b,L,S}\|}_2 \leq
    k_S{\|\muh^{a,L} - \muh^{b,L}\|}_2.
  \end{equation*}
  Therefore
  \begin{equation}\label{eq:proof_ineq}
    {\|\muh_i^{y + \varepsilon e_i,S} - \muh_i^{y - \varepsilon e_i,S}\|}_2 \leq
    k_S{\|\muh_i^{y + \varepsilon e_i,L} - \muh_i^{y - \varepsilon e_i,L}\|}_2,
  \end{equation}
  but since $\mathcal{M}_L$ constitutes an orthogonal linear projection of $y$ to $\muh^{y,L}$, there exist a $n \times n$ projection matrix (Hermitian and idempotent) $P_L$ such that
  \begin{equation*}
    \muh^{y,L} \equiv P_L y,
  \end{equation*}
and so we may further develop the right hand side of the
inequality~\eqref{eq:proof_ineq}
  \begin{align*}
    k_S{\|\muh_i^{y + \varepsilon e_i,L} - \muh_i^{y - \varepsilon e_i,L}\|}_2 =& k_S{\| P_L (y + \varepsilon e_i) - P_L (y - \varepsilon e_i) \|}_2
    = 2\varepsilon k_S {\| P_L e_i \|}_2 \\ =& 2\varepsilon k_S \sqrt{ e_i^T P_L^T P_L e_i } = 2\varepsilon k_S {P_L}_{ii}.
  \end{align*}
  On the other hand
  \begin{align*}
    \muh_i^{y + \varepsilon e_i,L} - \muh_i^{y - \varepsilon e_i,L} =& e_i^T P_L (y + \varepsilon e_i) - e_i^T P_L (y - \varepsilon e_i) \\
    =& 2\varepsilon e_i^T P_L e_i = 2\varepsilon {P_L}_{ii}.
  \end{align*}

  In summary, for every value of $\varepsilon \in \mathds{R}$, and for each $ i \in \{1,\dots,n\}$, we have shown that
  \begin{equation*}
    \muh_i^{y + \varepsilon e_i,S} - \muh_i^{y - \varepsilon e_i,S} \leq
    \muh_i^{y + \varepsilon e_i,L} - \muh_i^{y - \varepsilon e_i,L}.
  \end{equation*}
  This implies that every main-diagonal Jacobian component is smaller for the projection to $S$ than it is for the projection to $L$
  \begin{equation*}
    \frac{\partial{\muh_i^{y,S}}}{\partial{y_i}} \leq
    \frac{\partial{\muh_i^{y,L}}}{\partial{y_i}}\;,\; \mbox{ hence }
        \sum_{i=1}^n {\frac{\partial{\muh_i^{y,S}}}{\partial{y_i}}} \leq
    \sum_{i=1}^n {\frac{\partial{\muh_i^{y,L}}}{\partial{y_i}}}.
  \end{equation*}
  Because this is true for any observed data $y$, it is also true in expectation
  \begin{equation*}
    \E \Biggl( \sum_{i=1}^n {\frac{\partial{\muh_i^{y,S}}}{\partial{y_i}}} \Biggr) \leq
    \E \Biggl( \sum_{i=1}^n {\frac{\partial{\muh_i^{y,L}}}{\partial{y_i}}} \Biggr),
  \end{equation*}
  which, by Stein's lemma leads to $\omega_S \leq \omega_L$.
\end{proof}

\bibliographystyle{chicago}
\bibliography{optimism}

\begin{thebibliography}{}

\bibitem[\protect\citeauthoryear{Akaike}{Akaike}{1974}]{Akaike1974}
Akaike, H. (1974).
\newblock A new look at the statistical model identification.
\newblock {\em IEEE Transactions on Automatic Control\/}~{\em 19\/}(6),
  716--723.

\bibitem[\protect\citeauthoryear{Boyd and Vandenberghe}{Boyd and
  Vandenberghe}{2004}]{Boyd2009}
Boyd, S. and L.~Vandenberghe (2004).
\newblock {\em Convex optimization}.
\newblock Cambridge University Press.

\bibitem[\protect\citeauthoryear{Cherkassky and Mulier}{Cherkassky and
  Mulier}{2007}]{Cherkassky2007}
Cherkassky, V. and F.~Mulier (2007).
\newblock {\em Learning from data: Concepts, theory, and methods}.
\newblock Wiley-IEEE Press.

\bibitem[\protect\citeauthoryear{Davidov}{Davidov}{2006}]{Davidov2006}
Davidov, O. (2006).
\newblock Constrained estimation and the theorem of kuhn-tucker.
\newblock {\em Advances in Decision Sciences\/}~{\em 2006}.

\bibitem[\protect\citeauthoryear{Efron}{Efron}{1983}]{Efron1983}
Efron, B. (1983).
\newblock Estimating the error rate of a prediction rule: improvement on
  cross-validation.
\newblock {\em Journal of the American Statistical Association\/}~{\em
  78\/}(382), 316--331.

\bibitem[\protect\citeauthoryear{Efron}{Efron}{2004}]{Efron2004}
Efron, B. (2004).
\newblock The estimation of prediction error.
\newblock {\em Journal of the American Statistical Association\/}~{\em
  99\/}(467), 619--632.

\bibitem[\protect\citeauthoryear{Efron, Hastie, Johnstone, and
  Tibshirani}{Efron et~al.}{2004}]{Efron2004LARS}
Efron, B., T.~Hastie, I.~Johnstone, and R.~Tibshirani (2004).
\newblock Least angle regression.
\newblock {\em The Annals of statistics\/}~{\em 32\/}(2), 407--499.

\bibitem[\protect\citeauthoryear{Hastie, Tibshirani, and Friedman}{Hastie
  et~al.}{2009}]{Hastie2009}
Hastie, T., R.~Tibshirani, and J.~Friedman (2009).
\newblock {\em The elements of statistical learning: data mining, inference,
  and prediction}.
\newblock Springer Verlag.

\bibitem[\protect\citeauthoryear{Hoerl}{Hoerl}{1962}]{Hoerl1962}
Hoerl, A. (1962).
\newblock Application of ridge analysis to regression problems.
\newblock {\em Chemical Engineering Progress\/}~{\em 58\/}(3), 54--59.

\bibitem[\protect\citeauthoryear{Kato}{Kato}{2009}]{kato2009degrees}
Kato, K. (2009).
\newblock On the degrees of freedom in shrinkage estimation.
\newblock {\em Journal of Multivariate Analysis\/}~{\em 100\/}(7), 1338--1352.

\bibitem[\protect\citeauthoryear{Kattumannil}{Kattumannil}{2009}]{kumar2009stein}
Kattumannil, S. (2009).
\newblock On stein's identity and its applications.
\newblock {\em Statistics \& Probability Letters\/}~{\em 79\/}(12), 1444--1449.

\bibitem[\protect\citeauthoryear{Kr{\"a}mer and Sugiyama}{Kr{\"a}mer and
  Sugiyama}{2011}]{Kramer2011}
Kr{\"a}mer, N. and M.~Sugiyama (2011).
\newblock The degrees of freedom of partial least squares regression.
\newblock {\em Journal of the American Statistical Association\/}~{\em
  106\/}(494), 697--705.

\bibitem[\protect\citeauthoryear{Mallows}{Mallows}{1973}]{Mallows1973}
Mallows, C. (1973).
\newblock Some comments on cp.
\newblock {\em Technometrics\/}~{\em 15\/}(4), 661--675.

\bibitem[\protect\citeauthoryear{Schwarz}{Schwarz}{1978}]{Schwarz1978}
Schwarz, G. (1978).
\newblock Estimating the dimension of a model.
\newblock {\em The annals of statistics\/}~{\em 6\/}(2), 461--464.

\bibitem[\protect\citeauthoryear{Stein}{Stein}{1981}]{Stein1981}
Stein, C. (1981).
\newblock Estimation of the mean of a multivariate normal distribution.
\newblock {\em The Annals of Statistics\/}~{\em 9\/}(6), 1135--1151.

\bibitem[\protect\citeauthoryear{Stone}{Stone}{1974}]{Stone1974}
Stone, M. (1974).
\newblock Cross-validatory choice and assessment of statistical predictions.
\newblock {\em Journal of the Royal Statistical Society. Series B
  (Methodological)\/}, 111--147.

\bibitem[\protect\citeauthoryear{Tibshirani}{Tibshirani}{1996}]{Tibshirani1996}
Tibshirani, R. (1996).
\newblock Regression shrinkage and selection via the lasso.
\newblock {\em J. of the Royal Statistical Society: Series B\/}~{\em 58\/}(1),
  267--288.

\bibitem[\protect\citeauthoryear{Tibshirani and Taylor}{Tibshirani and
  Taylor}{2011}]{Tibshirani2011}
Tibshirani, R. and J.~Taylor (2011).
\newblock The solution path of the generalized lasso.
\newblock {\em The Annals of Statistics\/}~{\em 39\/}(3), 1335--1371.

\bibitem[\protect\citeauthoryear{Vapnik}{Vapnik}{2000}]{Vapnik2000}
Vapnik, V. (2000).
\newblock {\em The nature of statistical learning theory}.
\newblock Springer Verlag.

\bibitem[\protect\citeauthoryear{Wahba}{Wahba}{1983}]{Wahba1983}
Wahba, G. (1983).
\newblock Bayesian ``confidence intervals'' for the cross-validated smoothing
  spline.
\newblock {\em Journal of the Royal Statistical Society, Series B\/}~{\em
  45\/}(1), 133--150.

\bibitem[\protect\citeauthoryear{Ye}{Ye}{1998}]{Ye1998}
Ye, J. (1998).
\newblock On measuring and correcting the effects of data mining and model
  selection.
\newblock {\em Journal of the American Statistical Association\/}~{\em
  93\/}(441), 120--131.

\bibitem[\protect\citeauthoryear{Zou and Hastie}{Zou and
  Hastie}{2005}]{zou2005regularization}
Zou, H. and T.~Hastie (2005).
\newblock Regularization and variable selection via the elastic net.
\newblock {\em Journal of the Royal Statistical Society: Series B (Statistical
  Methodology)\/}~{\em 67\/}(2), 301--320.

\bibitem[\protect\citeauthoryear{Zou, Hastie, and Tibshirani}{Zou
  et~al.}{2007}]{Zou2007}
Zou, H., T.~Hastie, and R.~Tibshirani (2007).
\newblock On the degrees of freedom of the lasso.
\newblock {\em The Annals of Statistics\/}~{\em 35\/}(5), 2173--2192.

\end{thebibliography}

\centering


\section*{Supplementary material}

# When Does More Regularization Imply Fewer Degrees of Freedom? Sufficient Conditions and Counter Examples from Lasso and Ridge Regression

## Supplementary Material

Shachar Kaufman[*][1] and Saharon Rosset[†][1]

[1]Department of Statistics and Operations Research, Tel Aviv University

November 12, 2013

## S1 Use of optimism for practical model selection

The importance of Stein's lemma (6) is both theoretical and practical. Various authors have successfully used the lemma to find unbiased estimates for the optimism of popular models, such as shape restricted regression (Meyer and Woodroofe, 2000), the lasso (Zou et al., 2007) and support vector regression (Gunter and Zhu, 2007). Practical use of optimism theory in model selection is possible when one can either directly calculate $\omega$ or estimate it, e.g., using Eq. (6). In such cases, one only has to use (4) and an estimate of the expected training error such as the natural unbiased estimate $\hat{E}_{train} = \mathcal{E}(\hat{\mu}, y)$, to obtain an estimate of the expected prediction error: $\hat{E}_{pred} = \hat{E}_{train} + \hat{\omega}$. Specifically, this leads to Stein's unbiased risk estimate (SURE):

$$\hat{E}_{pred} = \frac{1}{n}\sum_{i=1}^{n} \left(\hat{\mu}_i - y_i\right)^2 + \frac{2\sigma^2}{n}\sum_{i=1}^{n} \frac{\partial \hat{\mu}_i}{\partial y_i}.$$

An unknown $\sigma$ is usually replaced by an estimate from the unregularized model.

This leads, for example to the following unbiased estimate for the prediction error of linear regression

$$\hat{E}_{pred} = \frac{1}{n}\sum_{i=1}^{n} \left(\hat{\mu}_i - y_i\right)^2 + \frac{2p\sigma^2}{n},$$

which coincides with Mallows' well-known $C_p$ model selection criterion (Mallows, 1973), where the second term is often viewed as a penalty on the degrees of freedom "used up" during fitting.

## S2 Toy example

To illustrate that the phenomenon of having more optimism in a nested modeling approach can translate into a significant effect on prediction performance, consider the following example.

---

[*]shachark@post.tau.ac.il
[†]saharon@post.tau.ac.il

Let $y \in \mathbb{R}^2$ be generated as $y_1 \sim Be(0.5)$ and $y_2 = 1$. Compare the Euclidean projections to the nested sets $S \subset L$, where $L$ is the triangle with vertices $(0,0), (0.5, 0.5), (0, 1)$, and $S$ is the triangle's bottom edge, i.e., the line segment from $(0,0)$ to $(0,1)$. The projection from $y$ to $L$ is fixed, $\hat{\mu} = (0.5, 0.5)$, thus the training error is 0.25 and the optimism is zero. On the other hand, the projection to $S$ gives $\hat{\mu} = (0, y_1)$ and has training error 0.5 and optimism $\text{cov}(y_1, \hat{\mu}_1) = \text{var}(y_1) = 0.25$ (see Figure S7, left panel for a graphical presentation).

## S3 "In the wild" example: regularized isotonic regression

While instructive and in many ways realistic, the counter examples in Section 4 of the main text are specifically constructed to show that significant violations of nondecreasing monotonicity between fitting and optimism can occur in nested setups, including in familiar ones. An important question is to what extent may we expect the "more fitting, less optimism" phenomenon to be prevalent "in the wild". In fact, we have originally encountered this phenomenon by coincidence while working on an algorithm for isotonic regression: Isotonic recursive partitioning.

*Isotonic regression (IR)* is a nonparametric modeling approach defined over a covariate space that has a notion of a *partial order*. A partial order is a reflexive, anti-symmetric, and transitive binary relation of points $x_i$ and $x_j$ in the covariate space, denoted $x_i \preceq x_j$, that indicates that for certain pairs of elements in the covariate space, one of the elements "precedes" the other. A primary example is the coordinate-wise order in $\mathbb{R}^m$, where points $x_i \in \mathbb{R}^m$ and $x_j \in \mathbb{R}^m$ are ordered if and only if they are ordered in every coordinate $x_{ik} \leq x_{jk}$, $k = 1, 2, \ldots, m \Leftrightarrow x_i \preceq x_j$. An *isotonic function* $f(x)$ is a function whose values are "correctly" ordered $f(x_i) \leq f(x_j)$ for every pair of partially ordered covariate points $x_i \preceq x_j$. Isotonic regression minimizes an error criterion, say squared loss (1), over the set of candidate fits $\tilde{\mu}$ that are isotonic functions of $x$. Formally:

$$\hat{\mu} = \underset{\tilde{\mu}}{\operatorname{argmin}} \Big\{ \|y - \tilde{\mu}\|_2^2, \text{ s.t. } \forall x_i \preceq x_j, \ \hat{\mu}_i \leq \hat{\mu}_j \Big\}. \tag{S3}$$

For a detailed review of isotonic regression, see Robertson et al. (1988).

Like other non-parametric approaches, IR models tend to overfit, especially in high dimension, and this calls for some form of regularization. Examining the structure of the IR problem, it can be shown that the solution is piecewise constant over a partition of the covariate space to a finite number $K$ of "hole-free" *blocks*. The recently proposed *Isotonic Recursive Partitioning* algorithm (IRP) (Luss et al., 2012) offers an efficient way of solving the isotonic regression problem based on an iterative partitioning scheme. There it is proven that a greedy approach, that starts with the entire covariate space as one block and iteratively splits it into blocks of decreasing size, is guaranteed to converge to the global solution of (S3) in $K$ iterations. It is proposed that by stopping at iteration $k < K$ along the *IRP path* of iterative splitting, IRP may be viewed as a family of nested modeling approaches with $k$ the tuning parameter (a larger $k$ necessarily leads to a finer fit, until $K$ is reached). This family is nested under Definition 2, since the IRP model with $k+1$ blocks is a refinement of the $k$-block model for any $k$.

Motivated by applications in genetics, Luss et al. (2012) were particularly interested in cases where the covariate space is a cartesian product of ternary-valued *genotypes* $x_i \in \{0, 1, 2\}^m$. For many traits such as height or risk of developing disease there is often reason to believe that the expected value of interest $\mathbb{E} y$ is an isotonic function of genotypes with the coordinate-wise partial order (Song and Nicolae, 2009). Here we show results of simulations with a simple three-dimensional ($m = 3$) sum-of-squares model:

$$y_i \sim \mathcal{N}\Big(\sum_{j=1}^{3} x_{ij}^2, \sigma^2\Big),$$

Figure S6 shows one such example with 256 observations and two choices for $\sigma$. Figure S7 (right panel) shows the relative contributions of training errors and optimism to prediction error for the setup showing the reversal in optimism. The effective degrees of freedom of models along the IRP path can be estimated empirically from repeated simulations directly or using the optimism theorem (5). Because of the true isotonic form of the generating function, with high probability the unconstrained fit, which fits the average to each of the $3^3 = 27$ possible genotypes, is itself isotonic. Consequently, the solution to the IR problem (S3) has *about* 27 effective degrees of freedom. The effective degrees of freedom of the full IR fit (which is always the last fit on the IRP path) can be unbiasedly estimated by the number of blocks in its associated partition (see Meyer and Woodroofe, 2000). Because of the true isotonic form of the generating function, here with high probability the number of blocks is the size of the finite covariate space $3^3 = 27$. All effective degrees of freedom profiles arrive eventually at this value, but while we expected to see the regularization reflected as a monotonically nondecreasing profile, for many of the models examined this did not occur – we consistently obtained optimism profiles similar to those of Example 4.2, where regularized models had significantly higher effective degrees of freedom.

The reason for the extreme violation of monotonicity of optimism in nesting evident in this example is the highly non-convex "star shape" nature of the model sets in use here: *all models with up to k isotonic blocks*. The geometry of these sets makes it that much easier to encounter true data distributions that are situated similarly to examples 2.2, 4.1 and 4.2, and the same mechanisms are behind the erratic behavior of the optimism profiles.

## S4 Proof of Proposition 1

Suppose observations $y_i$ which are components of the vector $y_{n \times 1}$ are mutually uncorrelated but not homoscedastic

$$\Lambda_{ij} = \begin{cases} \sigma_i^2 & \text{, if } i = j \\ 0 & \text{, otherwise} \end{cases}$$

For a general linear smoother model where $S$ is any $n \times n$ matrix uncorrelated with the data $y$, the expected optimism is

$$\omega = \frac{2}{n}\sum_{i=1}^n \text{cov}(y_i, \hat{\mu}_i) = \frac{2}{n}\sum_{i=1}^n \text{cov}(y_i, (Sy)_i) = \frac{2}{n}\text{tr}(\text{cov}(y, Sy)) = \frac{2}{n}\text{tr}(\Lambda S)$$

Specifically for penalized ridge regression $S$ is given by Eq. (7). With the SVD $X = VDU^T$ we have

$$\text{tr}(\Lambda S) = \text{tr}\big(\Lambda U D V^T (V D U^T U D V^T + \lambda I)^{-1} V D U^T\big) =$$
$$\text{tr}(\Lambda U D V^T (V D^2 V^T + \lambda V V^T)^{-1} V D U^T) =$$
$$\text{tr}\Big(\Lambda U \, \text{diag}\Big(\frac{d_j^2}{d_j^2 + \lambda}\Big) U^T\Big) = \sum_{i=1}^n \sigma_i^2 \sum_{j=1}^p \frac{d_j^2}{d_j^2 + \lambda} u_{ij}^2$$

Meaning the expected optimism still shrinks with growing $\lambda$ as in the homoscedastic case.

## References

Gunter, L. and J. Zhu (2007). Efficient computation and model selection for the support vector regression. *Neural Computation* 19(6), 1633–1655.

Luss, R., S. Rosset, and M. Shahar (2012). Efficient regularized isotonic regression with application to gene–gene interaction search. *The Annals of Applied Statistics 6*(1), 253–283.

Mallows, C. (1973). Some comments on cp. *Technometrics 15*(4), 661–675.

Meyer, M. and M. Woodroofe (2000). On the degrees of freedom in shape-restricted regression. *The Annals of Statistics 28*(4), 1083–1104.

Robertson, T., F. Wright, R. Dykstra, and T. Robertson (1988). *Order restricted statistical inference*, Volume 229. Wiley New York.

Song, M. and D. Nicolae (2009). Restricted parameter space models for testing gene-gene interaction. *Genetic epidemiology 33*(5), 386–393.

Zou, H., T. Hastie, and R. Tibshirani (2007). On the degrees of freedom of the lasso. *The Annals of Statistics 35*(5), 2173–2192.

# S5 Supplementary figures

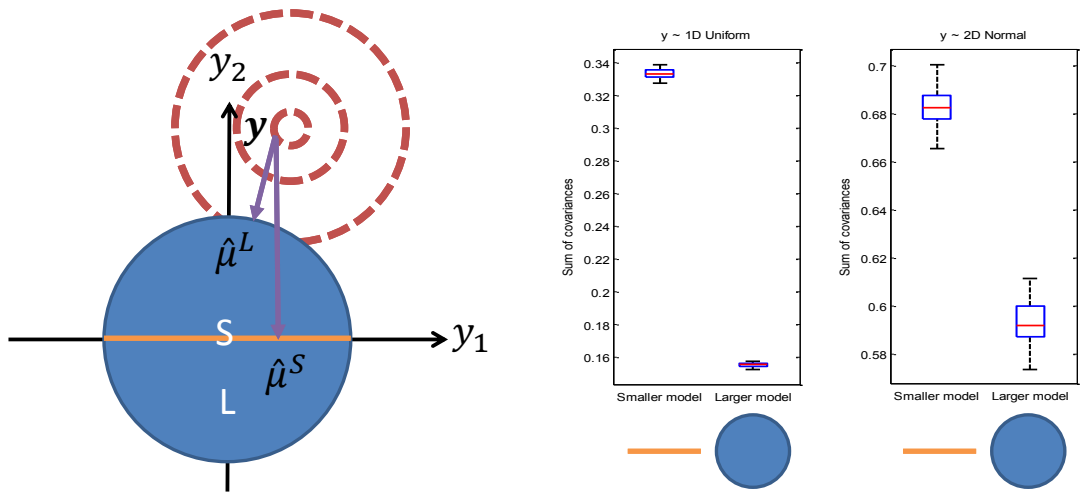

Figure S1: Toy example. Left: The geometrical setup, presented with an i.i.d. normal (dashed) data distribution. Right: A comparison of the optimism of the two modeling approaches projecting onto $S$ and $L$.

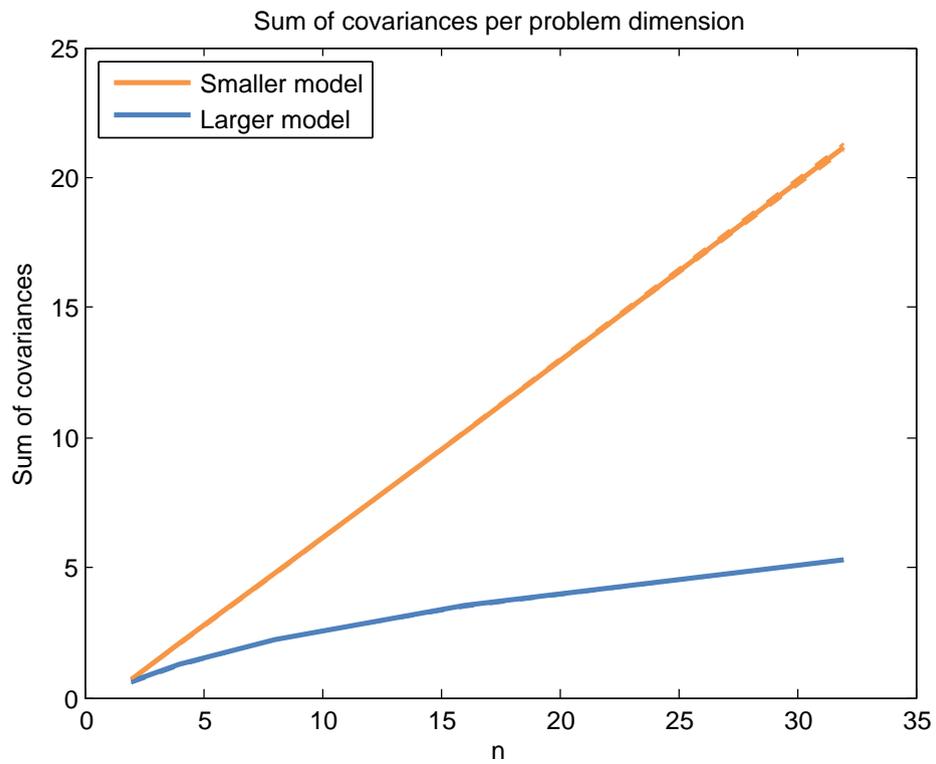

Figure S2: Evolution of the optimism as the dimension grows in the toy example of Figure S1 (with normally distributed data). Shown are the mean (solid) and 95% confidence intervals (dashed, negligible) for the sum of covariances from Eq. (5).

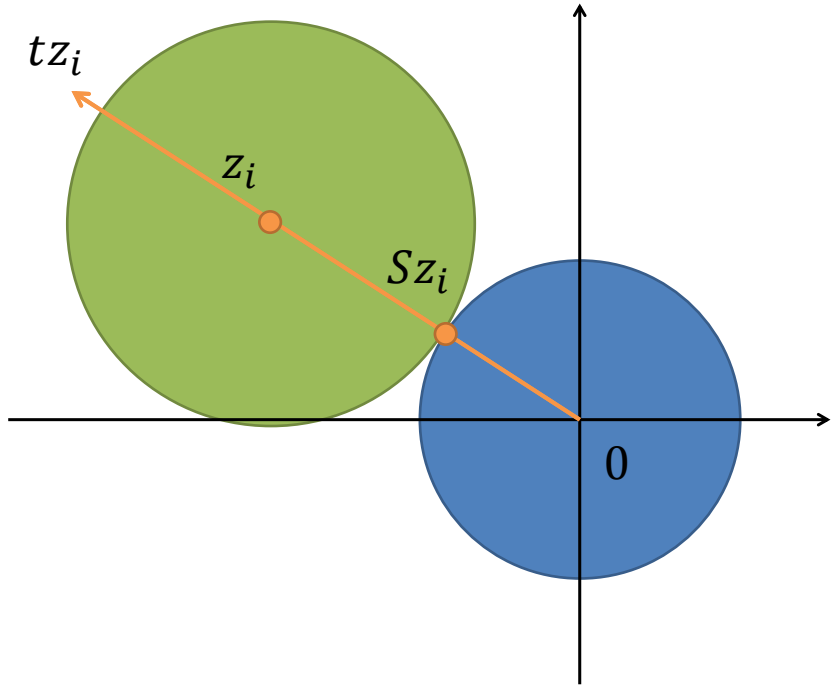

Figure S3: If $\hat{\mu}^S(y) = Sy$ satisfies Definition 2 over a convex set: (i) the convex set has to include the origin, (ii) all eigenvalues of $S$ must be in $[0, 1]$, (iii) and since $\hat{\mu}^L = Ly$ is a projection on a larger set, $Lz_i$ must be closer to $z_i$ than $Sz_i$ is, thus $\|Lz_i\|_2^2 \geq \|Sz_i\|_2^2$.

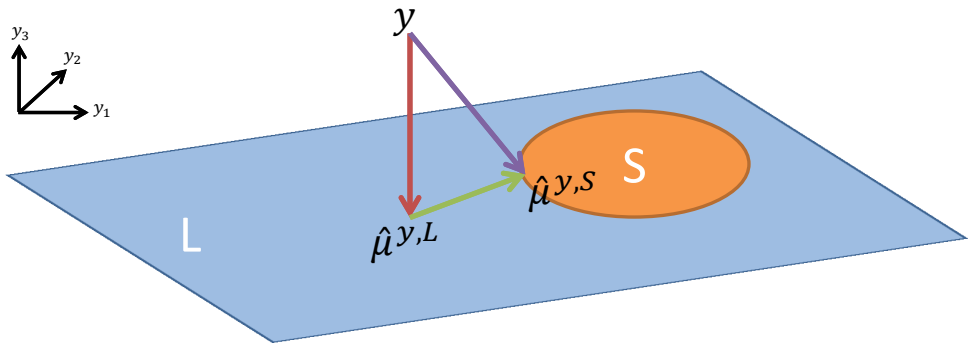

Figure S4: When $L$ is a linear space, projecting a point to $S \subseteq L$ may be achieved by first linearly projecting onto $L$ and from that point projecting onto $S$.

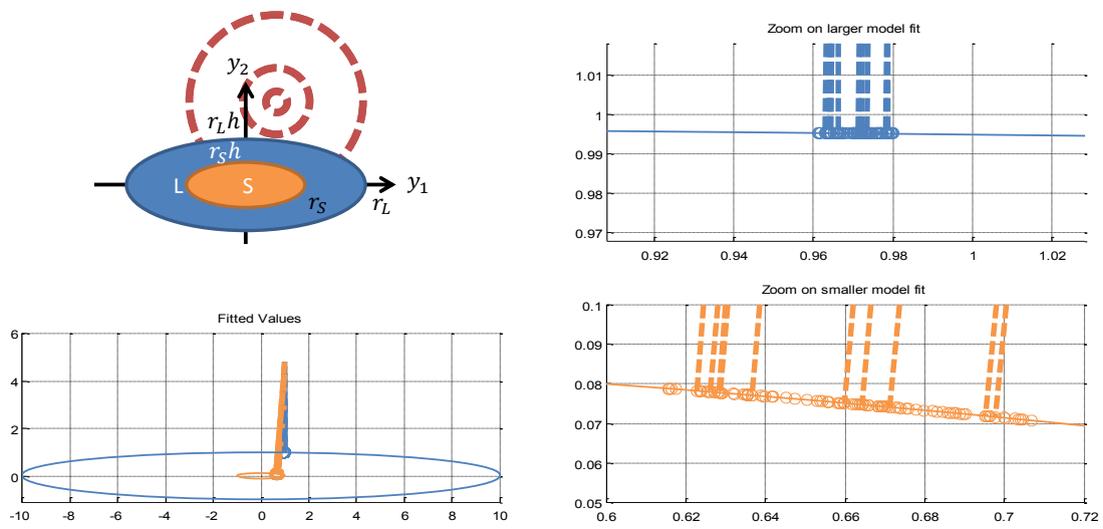

Figure S5: Constrained ridge regression counterexample, setup (top left) and simulation results.

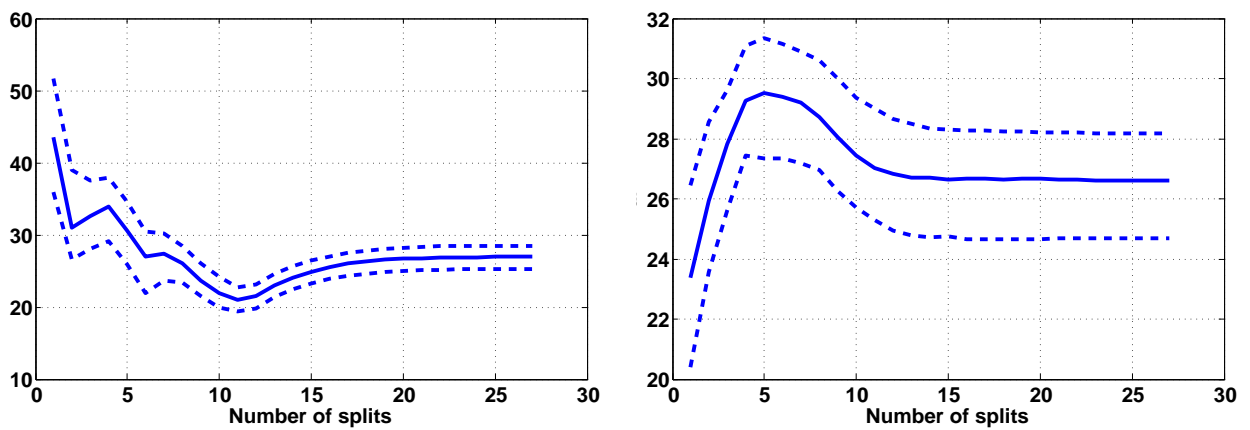

Figure S6: Isotonic recursive partitioning counterexample. The effective degrees of freedom profile for the entire regularization path is shown. Solid curves depict the mean, dashed curves give 95% confidence intervals. Left panel: $\sigma^2 = 1$; Right panel: $\sigma^2 = 10$.

8

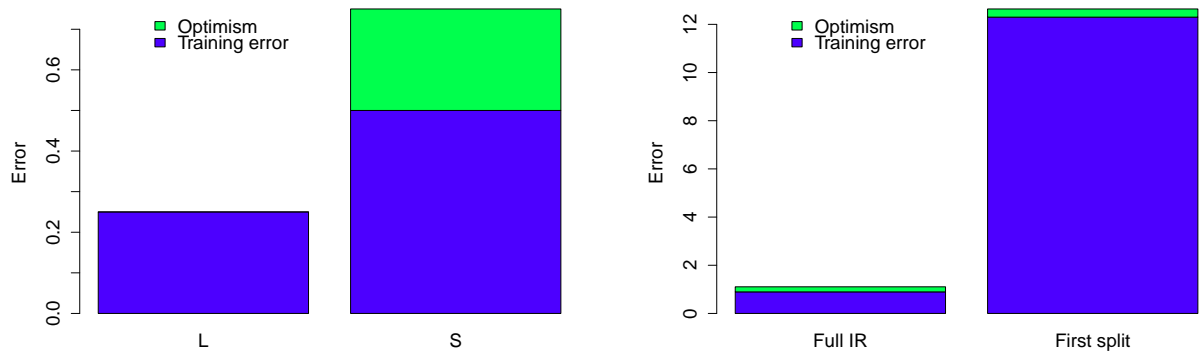

Figure S7: Training error, optimism, and prediction error. Left: toy example of Section S2. Right: IRP example matching the left panel of Figure S6.



\end{document}